\documentclass[11pt]{amsart}
\usepackage{graphicx}
\allowdisplaybreaks
\usepackage{amsmath,amssymb,mathrsfs,amsthm}

\usepackage{verbatim}
\usepackage[spanish,USenglish]{babel} 
\usepackage{color}
\usepackage{tikz}
\usepackage{graphicx}

\newtheorem{thm}{Theorem}[section]
\newtheorem{corollary}[thm]{Corollary}
\newtheorem{proposition}[thm]{Proposition}
\newtheorem{theorem}[thm]{Theorem}
\newtheorem{lemma}[thm]{Lemma}

\theoremstyle{definition}

\theoremstyle{remark}
\newtheorem{remark}[thm]{Remark}

\numberwithin{equation}{section}

\def\N{\mathbb{N}}
\def\R{\mathbb{R}}

\def\C{\mathbb{C}}

\def\Z{\mathbb{Z}}

\def\T{\mathbb{T}}

\def\Tn2{{\mathcal T}^{(n)}_2}

\def\R{\mathbb{R}}
\def\N{\mathbb{N}}
\def\C{\mathbb{C}}

\begin{document}

\title[Catalan generating functions for bounded operators]
{Catalan generating functions for bounded operators}

\author[Miana]{Pedro J. Miana}
\address{Departamento de Matem\'aticas, Instituto Universitario de Matem\'aticas y Aplicaciones, Universidad de Zaragoza, 50009 Zaragoza, Spain.}
\email{pjmiana@unizar.es}

\author[Romero]{Natalia Romero}
\address{Departamento de Matem\'aticas y Computaci\'on,  Universidad de la Rioja, 26006 Logro\~{n}o, Spain.}
\email{natalia.romero@unirioja.es}

\thanks{ Pedro J. Miana has been partially supported by Project ID2019-105979GBI00, DGI-FEDER, of the MCEI and Project E48-20R, Gobierno de Arag\'on, Spain. Natalia Romero   has been partially supported by the project MTM2018-095896-B-C21
of the Spanish Ministry of Science.}



\maketitle

\begin{abstract}In this paper we study the solution of the quadratic equation $TY^2-Y+I=0$ where $T$ is a linear and bounded operator on a Banach space $X$. We describe the spectrum set and the resolvent operator of $Y$ in terms of operator $T$.
In the case that $ 4T$ is a power-bounded operator, we show that a solution (named Catalan generating function) is given by the Taylor series
$$
C(T):=\sum_{n=0}^\infty C_nT^n,
$$
where the sequence $(C_n)_{n\ge 0}$ is the well-known Catalan numbers. We express $C(T)$ by means of an integral representation which involves the resolvent operator $(\lambda-T)^{-1}$. Some particular examples to illustrate our results are given, in particular an iterative method defined  for square matrices $T$ which involves Catalan numbers. \end{abstract}

\section{Introduction}
The well-known Catalan numbers $(C_n)_{n\ge 0}$  given by the formula
$$
 C_n={1\over n+1}{2n\choose n},\quad  \ n\ge 0,
$$
appear in a wide range of problems. For instance, the Catalan number $C_n$ counts the  number of ways to triangulate a regular polygon with $n +2$ sides; or,
the number of  ways  that   $2n$ people seat around a circular table are simultaneously
shaking hands with another person at the table in such a way that none of the arms cross each other, see for example \cite{[Sl1], [St2]}.
They have been studied in depth in many papers and monographs (see for example \cite{[CC],[MR2], [Sh74], [St2]}) and the Catalan
sequence is probably the most frequently encountered sequence.

The generating function of the Catalan sequence $c=(C_n)_{n\ge 0}$ is defined by
\begin{equation}\label{gene}
C(z):=\sum_{n=0}^\infty C_nz^n= {1- \sqrt{1-4z}\over 2z}, \quad z\in D(0,{1\over 4}):=\{z\in \C\,\,|\,\, \vert z\vert<{1\over 4}\}.
\end{equation}
This function satisfies the quadratic equation
$
zy^2-y+1=0.
$ The main object of this paper is to consider this quadratic equation in the set of linear and bounded operators, ${\mathcal B}(X)$ on a Banach space $X$, i.e.,
\begin{equation}\label{Ceq}
TY^2-Y+I=0,
\end{equation}
where $I$ is the identity on the Banach space, and  $T , Y \in {\mathcal B}(X)$.  Formally, some solutions of this vector-valued quadratic equations are expressed by
$$
Y={1\pm\sqrt{1-4T}\over 2T},
$$
 which involves the (non-trivial) problems of the square root of operator $1-4T$ and the inverse of operator $T$.

  In general, the equation (\ref{Ceq}) may have no solution, one, several or infinite solutions, see examples in Section 6.  Note that study of quadratic equations in Banach space $X$ which $\hbox{dim}(X)\ge 2$ is much complicated than in the scalar case. For example
 there are infinite symmetric square roots of $I_2 \in \R^{2\times 2}$ given by
$$
\begin{pmatrix}
1 & 0\\
 0 & 1
\end{pmatrix}, \qquad
\begin{pmatrix}
-1 & 0\\
 0 & -1
\end{pmatrix}, \qquad
 \begin{pmatrix}
a & \pm\sqrt{1-a^2}\\
 \pm\sqrt{1-a^2} & -a
\end{pmatrix}
$$
with $a\in[-1,1]$.

As far as we are aware, no useful necessary and sufficient conditions
for the existence of solution of quadratic equations in Banach spaces are known, even in the classical
case of square roots in finite-dimensional spaces. To find some
easily applicable conditions is of interest, in part because these equations
are frequently used in the study of, for example, physical or biological phenomena.

In 1952, Newton's method was generalized to Banach
space by Kantorovich. Kantorovich's theorem asserts that the iterative method of Newton, applied to a most general system of nonlinear equations $P(x) =0$, converges to a solution $x^*$ near some given point $x_0$, provided the Jacobian of the system satisfies a Lipschitz condition near $x_0$  and its inverse at $x_0$ satisfies certain boundedness conditions. The theorem also gives computable error bounds for the iterates. From here, a large theory has been developed to obtain sharp iterative methods to approximate solutions of non linear equation (see for example \cite{HRR, LAN, McF}) and in particular quadratic  matrix equations (\cite{HR, LAN85}).

The paper is organized as follows. In the second section, we  show new results about the well-known Catalan numbers sequence $(C_n)_{n\ge 0}$. In Theorem  \ref{integrals}, we prove that following technical identity holds
$$
\int_0^\infty{\sqrt{t}\over (t+1)(t+z)^{j+1}}dt={\pi\over 2\sqrt{z}(z-1)^{j}}\sum_{k=j}^{\infty}C_k\left({z-1\over 4z}\right)^k, \qquad \Re(z)\ge {1\over 2},
$$
for $j\ge 1$. A nice result about solutions of quadratic equations
$$
xy^2-y+1=0, \qquad -xz^2-z+1=0,
$$
are given in Theorem \ref{quartic}: the arithmetic mean $\displaystyle{y+z\over 2}$ is solution of the  biquadratic equation $4x^2w^4-w^2+1=0$.

We consider the sequence $c=(C_n)_{n\ge 0}$ as an element in the Banach algebra $\ell^1(\N^0, {1\over 4^n})$ in the third section. We describe the spectrum set $\sigma(c)$ in  Proposition \ref{spec} and the resolvent element  $(\lambda-c)^{-1}$ in  Theorem \ref{inver}.

In the forth section, we study spectral properties of the solution of quadratic equation (\ref{Ceq}) with $T\in {\mathcal B}(X)$. We prove several results between $\sigma (T)$ and $\sigma(Y)$ where $\sigma(\cdot)$ denotes the spectrum set of the operator $T$. Moreover, we express  $(\lambda-Y)^{-1}$ in terms of the resolvent of operator $T$ in  Theorem \ref{spectral}.

For operators $T$ which $4T$ are power-bounded, we define the generating Catalan function
$$
C(T):=\sum_{n\ge 0}C_nT^n.
$$
This operator solves the quadratic equation (\ref{Ceq}) and has interesting properties connected with $T$, see Theorem \ref{genne}; in particular the following integral representation holds,
$$
C(T)x={1\over \pi}\int_{1\over 4}^\infty {\sqrt{\lambda-{1\over 4}}\over \lambda}(\lambda-T)^{-1}x d\lambda, \qquad x\in X.
$$

 In the last section we illustrate our results with some  examples of operators $T$ in  the equation  (\ref{Ceq}). We consider the Euclidean space $\C^2$ and matrices $$T=\lambda I_2, \,\,  \begin{pmatrix}
0 & \lambda\\
\lambda & 0
\end{pmatrix}, \,\, \begin{pmatrix}
0 & \lambda\\
0 & 0
\end{pmatrix}.$$
We solve the equation (\ref{Ceq}) and  calculate $C(T)$ for these matrices. We also check $C(a)$ for some particular values of $a\in \ell^1(\N^0, {1\over 4^n})$ . Finally we present an iterative method for matrices $\R^{n\times n }$ which are defined using Catalan numbers.

\section{ Some news results about Catalan numbers}\label{GCHO}
The Catalan numbers may be defined recursively by $C_0=1$ and
\begin{equation}\label{cata}C_n= \sum_{i=0}^{n-1} C_i C_{n-1-i}, \qquad n\ge 1, \end{equation}
 and first terms in this sequence are
$
1,\,\,1,\,\,2,\,\,5,\,\,14, \,\,42, \,\,132, \dots .
$
The generating function of the Catalan sequence $c=(C_n)_{n\ge 0}$ is given in
(\ref{gene}). This function satisfies the quadratic equation
\begin{equation}\label{quadratic}
zC^2(z)-C(z)+1=0, \qquad z\in D(0,{1\over 4}),
\end{equation}
see for example \cite[Section 1.3]{[St2]}. The second solution of this quadratic equation  is given by
$$
{1\over zC(z)}={1+ \sqrt{1-4z}\over 2z}, \qquad z\in D(0,{1\over 4})\backslash\{0\}.
$$

The following theorem shows that  the arithmetic mean of two solutions of these quadratic equations is also solution of a  biquadratic equation, closer to the previous ones.

\begin{theorem}\label{quartic} Let $A$ be a commutative algebra over $\R$ or $\C$ with $x\in A$. If   $y$ and $z$ are solutions of the quadratic equations
$$
xy^2-y+1=0, \qquad -xz^2-z+1=0,
$$
then $\displaystyle{y+z\over 2}$ is a solution of the biquadratic equation $4x^2w^4-w^2+1=0$.

\end{theorem}
\begin{proof}Note that it is enough to show that $x^2(y+z)^4-(y+z)^2+4=0$. We write for a while
\begin{eqnarray*}
P:&=& (xy^2-y+1)(xy^2+4xyx+3xz^2)\cr
Q:&=&(-xz^2-z+1)(-xz^2-4xyx-3xz^2)\cr
P+Q&=&x^2(y+z)^4-xy^3-xy^2z-2xy^2+xz^3+xz^2y+2xz^2\cr
&=&x^2(y+z)^4-xy^2(y+z+2)+xz^2(z+y+2)\cr
&=&x^2(y+z)^4+(y+z+2)(xz^2- xy^2).
\end{eqnarray*}
Since $y$ and $z$ are solutions of these quadratic equations, we have that $xz^2- xy^2=2-(y+z)$ and
$$
0=P+Q= x^2(y+z)^4+(y+z+2)(2-(y+z))=x^2(y+z)^4-(y+z)^2+4.
$$
\end{proof}

\begin{remark} The sum of both equations $xy^2-y+1=0$ and $ -xz^2-z+1=0$ gives
$$
{y+z\over 2}-1={1\over 2} x(y^2-z^2)=2x\left({y-z\over 2}\right)\left({y+z\over 2}\right),
$$
and we may obtain  $\displaystyle{y-z\over 2}$ in terms of $\displaystyle{y+z\over 2}$  whatever the inverse of $\displaystyle{2x {y+z\over 2}}$ exists in the algebra $A$, i.e.
$$
{y-z\over 2}=\left({y+z\over 2}-1\right) \left(2x {y+z\over 2}\right)^{-1}  .
$$
\end{remark}

\medskip
As a direct aplication of Abel's theorem to  (\ref{gene}), we obtain that
\begin{equation}\label{norm}
\sum_{n=0}^\infty {C_n\over 4^n}=\lim_{z\to {1\over 4}}C(z)=2, \qquad \sum_{n=0}^\infty {C_n\over (-4)^n}=\lim_{z\to -{1\over 4}}C(z)=2(\sqrt{2}-1),
\end{equation}
(\cite[Exercise A.66]{[St2]}). In fact one has that
$$
C_n\sim{4^n\over \sqrt{\pi}n^{3\over 2}}, \qquad n\to \infty,
$$
(\cite[Exercise A.64]{[St2]}).

A straightforward consequence of the generating formula (\ref{gene}) and Theorem \ref{quartic} is the following proposition, where we consider the odd and even parts, $C_o(z),$ and $C_e(z)$ of function $C(z)$. The proof is left to the reader.

\begin{proposition}\label{descopo} Let $c=(C_n)_{n\ge 0}$ be the Catalan sequence. Then
\begin{eqnarray*}
C_e(z):=\sum_{n=0}^\infty C_{2n}z^{2n}&=&{\sqrt{1+4z}-\sqrt{1-4z}\over 4z},\cr
C_o(z):=\sum_{n=0}^\infty C_{2n+1}z^{2n+1}&=&{2-\sqrt{1+4z}-\sqrt{1-4z}\over 4z},\cr
\end{eqnarray*}
for $\vert z\vert \le {1\over 4}$. In particular, $4z^2C_e^4(z)-C_e(z)^2+1=0$, $C_e^2(z)=C(4z^2)$, $\displaystyle{C_o(z)}={C_e(z)-1\over 2zC_e(z)}$ for $\vert z\vert \le {1\over 4}$ and
$$
\sum_{n=0}^\infty {C_{2n}\over 4^{2n}}=\sqrt{2},\qquad
\sum_{n=0}^\infty {C_{2n+1}\over 4^{2n+1}}=2-\sqrt{2}.
$$
\end{proposition}

\medskip
Catalan numbers have several integral representations, for example
$$
C_n={1\over 2\pi}\int_0^4t^n{\sqrt{4-t\over t}}dt={2^{2n+1}\over \pi}\beta({3\over 2}, n+{1\over 2}),
$$
where the function $\beta $ is the well-known Euler Beta function, $\beta(u,v):=\int_0^1t^{u-1}(1-t)^{v-1}dt$ for $u,v>0$, see the monography \cite{[St2]} and the survey \cite{Qi}.  In the next theorem, we present a new results which involves the Taylor polynomials of the Catalan generating function $C(z)$.

\begin{theorem}\label{integrals} Given $1\not =z\in \C^{+}$, then
\begin{eqnarray*}
\int_0^\infty{\sqrt{t}\over (t+1)(t+z)}dt&=&{\pi\over z-1}\left({\sqrt{z}}-{1}\right),\cr
\int_0^\infty{\sqrt{t}\over (t+1)(t+z)^{j+1}}dt&=&{\pi\over (z-1)^{j+1}}\left({\sqrt{z}}-{1}-{z-1\over 2\sqrt{z}}\sum_{k=0}^{j-1}C_k\left({z-1\over 4z}\right)^k\right)\cr
&=&{\pi\over 2\sqrt{z}(z-1)^{j}}\sum_{k=j}^{\infty}C_k\left({z-1\over 4z}\right)^k,
\end{eqnarray*}
for $j\ge 1$ and where the last equality holds for $\Re(z)\ge {1\over 2}$.
\end{theorem}

\begin{proof} The first integral is a easy exercise of elemental calculus. To do the second one, note that
$$
{\sqrt{t}\over (t+1)(t+z)^{j+1}}={1\over z-1}\left({\sqrt{t}\over (t+1)(t+z)^{j}}- {\sqrt{t}\over (t+z)^{j+1}}\right),
$$ and then
\begin{eqnarray*}
\int_0^\infty{\sqrt{t}\over (t+1)(t+z)^{j+1}}dt&=&{1\over z-1}\left(\int_0^\infty{\sqrt{t}\over (t+1)(t+z)^{j}}dt- {\beta({3\over 2}, j-{1\over 2})\over z^{j-{1\over 2}}}\right)\cr
&=&{1\over z-1}\left(\int_0^\infty{\sqrt{t}\over (t+1)(t+z)^{j}}dt- {\pi C_{j-1}\over 2\sqrt{z}(4z)^{j-1}}\right),
\end{eqnarray*}
for $j\ge 1$. We iterate this formula to get the final expression.
\end{proof}

\begin{remark} By holomorphic property, Theorem \ref{integrals} holds for $z\in \C\backslash (-\infty, 0]$. Moreover for $z=1,$
$$
\int_0^\infty{\sqrt{t}\over (t+1)^{j+2}}={\pi C_{j}\over 2^{2j+1}}=\pi\lim_{z\to 1}\left({\sqrt{z}-{1}-{z-1\over 2\sqrt{z}}\sum_{k=0}^{j-1}C_k\left({z-1\over 4z}\right)^k\over  (z-1)^{j+1}}\right),
$$
for $j\ge 1$. Finally, when $j\to \infty$, we recover the generating formula
$$
\sum_{k=0}^\infty C_k\left({z-1\over 4z}\right)^k= {2\over \sqrt{z}+1}, \qquad \vert {{z-1}\over z}\vert \le 1.
$$
\end{remark}

\section{The sequence of Catalan numbers}\label{GCHO}
We may interpret the equality (\ref{norm}) in terms of norm in the weight Banach algebra $\ell^1(\N^0, {1\over 4^n}).$ This algebra is formed by sequence $a=(a_n)_{n\ge 0}$ such that
$$
\Vert a\Vert_{1, {1\over 4^n}}:=\sum_{n=0}^\infty {\vert a_n\vert\over 4^n}<\infty,
$$
and the product is the usual convolution $\ast $ defined by
$$
(a\ast b)_n=\sum_{j=0}^na_{n-j}b_j, \qquad a,b \in \ell^1(\N^0, {1\over 4^n}).
$$
The canonical base $\{\delta_j\}_{j\ge 0}$ is formed by sequences such that $(\delta_j)_n=\delta_{j,n}$ is the known delta Kronecker. Note that $\delta_1^{\ast n}= \delta_1\dots^{n}\delta_1=\delta_n$ for $n\in \N$.
This Banach algebra has identity element, $\delta_0$, its spectrum set is the closed disc $\overline{D(0,{1\over 4})}$  and its Gelfand transform is given by the $Z$-transform
$$
Z(a)(z):= \sum_{n=0}^\infty {a_n}z^n, \qquad z\in \overline{D(0,{1\over 4})}.
$$
It is straightforward to check that $Z(\delta_n)(z)=z^n$ for $n\ge 0$ (see, for example, \cite{La}).

In the next proposition, we collect some properties of the Catalan sequence $c$  in the language of the Banach algebra $\ell^1(\N^0, {1\over 4^n})$. In particular the identity (\ref{cata}) is equivalent to the item (iii).
\begin{proposition}\label{first} Take  $c=(C_n)_{n\ge 0}$. Then
\begin{itemize}
\item[(i)] $\Vert c\Vert_{1, {1\over 4^n}}=2$.
\item[(ii)]$C(z)=Z(c)(z)$ for   $z\in D(0,{1\over 4})$.
\item[(iii)] $\delta_1* c^{*2}-c+\delta_0=0$.
\end{itemize}

\end{proposition}

We recall that the  resolvent set of $a\in \ell^1(\N^*, {1\over 4^n})$, denoted as $\rho(a)$, is defined by $$\rho(a):=\{\lambda \in \C\, \, : \,\, (\lambda\delta_0-a)^{-1}\in \ell^1(\N^0, {1\over 4^n})\},
$$
and the spectrum set of  $a$  is denoted by $\sigma(a)$ and given by $\sigma(a):=\C\backslash \rho(a)$.

\medskip

\begin{proposition} \label{spec} The spectrum of the Catalan sequence $c=(C_n)_{n\ge 0}$ is given by
$
\sigma(c)=C(\overline{D(0,{1\over 4})})
$
and its boundary by $$\partial (\sigma(c))=\left\{2e^{-i\theta}\left(1-\sqrt{2\vert\sin({\theta\over 2})}\vert e^{i(\pi-\theta)\over 4}\right)\,\, :\,\,\theta\in (-\pi, \pi)\right\}.$$

\end{proposition}

\begin{proof}

As the algebra $\ell^1(\N^0, {1\over 4^n})$   has identity, we apply \cite[Theorem 3.4.1]{La} to get the equality set $\sigma(c)=C(\overline{D(0,{1\over 4})})$.

We write $\T:= \{e^{i\theta}\,\,:\,\, \theta\in (-\pi, \pi)\}$. Take $z\in
\partial (\sigma(c))= C({1\over 4}\T)
$
and
$$
z=2e^{-i\theta} \left(1-\sqrt{2i \sin ({\theta\over 2})e^{i\theta\over 2}}\right)= 2e^{-i\theta}\left(1-\sqrt{2\vert\sin({\theta\over 2})}\vert e^{i(\pi-\theta)\over 4}\right),
$$
for $\theta \in (-\pi, \pi).$
\end{proof}

\begin{remark} In the Figure 1, we plot the set $\partial (\sigma(c)).$

\begin{figure}
 \centering
   \includegraphics[scale=0.5]{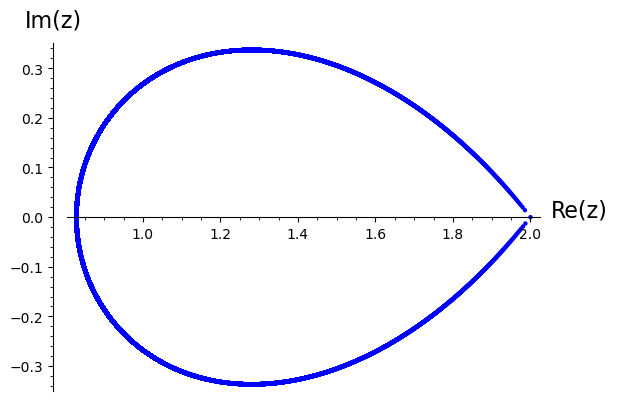}
   \caption{The set $\partial(\sigma(c))$.}
\end{figure}

\end{remark}

Given $\lambda \in \C$, we consider the geometric progression $p_{\lambda}:=´({1 \over \lambda^n})_{n\ge 0}.$ Note that $p_{\lambda} \in \ell^1(\N^0, {1\over 4^n})$ if and only if $\vert \lambda \vert >{1\over 4}$. Moreover
$$
(\lambda-\delta_1)^{-1}={1\over \lambda}p_{\lambda}, \qquad \vert \lambda \vert >  {1\over 4},
$$
and $Z((\lambda-\delta_1)^{-1})(z)=(\lambda-z)^{-1}$ for $z \in D(0,{1\over 4})$ and  $\vert\lambda \vert >  {1\over 4}.$ Note that
\begin{equation}\label{integra}
C_n={1\over \pi}\int_{1\over 4}^\infty{\sqrt{\lambda -{1\over 4}}\over \lambda^{2+n}}d\lambda={1\over \pi}\int_{1\over 4}^\infty{\sqrt{\lambda -{1\over 4}}\over \lambda^{2}}p_\lambda(n)d\lambda, \qquad n \ge 0,
\end{equation}
see, for example \cite[Section 4.7]{Qi}.

In the next theorem, we express $(\lambda-c)^{-1}$ in terms of $p_\lambda$ and $c$ for $\lambda \in\Omega$ where $\Omega$ is the following open set,
$$
\Omega:=\{\lambda \in \C\, \,:\,\,\left| { \lambda-1\over \lambda^2}\right|>{1\over 4}\}.
$$

\begin{theorem}\label{inver}  The inverse of the Catalan sequence $c$ is given  $c^{-1}= \delta_0-\delta_1\ast c$ and
$$
(\lambda-c)^{-1}={\delta_0\over \lambda}+{1\over \lambda(\lambda-1)}p_{\lambda-1\over \lambda^2}+{1\over \lambda^2 } c-{1\over \lambda^2}c\ast p_{\lambda-1\over \lambda^2}, \qquad \lambda \in \Omega\backslash\{0\}.
$$

\end{theorem}
\begin{proof} By Proposition \ref{first} (iii), $(\delta_0- \delta_1\ast c)\ast c= \delta_0$ and we conclude that $c^{-1}= \delta_0-\delta_1\ast c$ .

For $\lambda\in \Omega$,  we apply the Zeta transform to get that
$$
Z({\delta_0\over \lambda}+{1\over \lambda(\lambda-1)}p_{\lambda-1\over \lambda^2})(z)={1\over \lambda}\left(1+{1\over \lambda-1-z\lambda^2}\right)={z\lambda-1\over z\lambda^2-\lambda+1},
$$
$$
Z({1\over \lambda^2 } c-{1\over \lambda^2}c\ast p_{\lambda-1\over \lambda^2})(z)={C(z)\over \lambda^2}\left(1-{\lambda-1\over \lambda-1-z\lambda^2}\right)={zC(z)\over z\lambda^2-\lambda+1},
$$
for $z\in D(0,{1\over 4})$. To conclude the equality, we check that
$$
(\lambda-C(z))\left({z\lambda-1+zC(z)\over z\lambda^2-\lambda+1}\right)= {z\lambda^2-\lambda+C(z)-zC^2(z)\over z\lambda^2-\lambda+1}=1,
$$
where we have applied the quadratic identity (\ref{quadratic}).
\end{proof}

\begin{remark} As Theorem \ref{inver} shows, the set $\Omega$ is strictly contained in $\rho(c)$. Moreover the boundary of $\sigma(c)$ is contained in the boundary of $\Omega$, i.e. $\partial(\sigma(c))\subset \partial(\Omega)=\{\lambda \in \C\, \,:\,\,\left| { \lambda-1\over \lambda^2}\right|={1\over 4}\}$. In the Figure 2, we plot both sets, $\partial(\Omega)$ in blue and
$\partial(\sigma(c))$ in red.

\begin{figure}
  \centering
   \includegraphics[scale=0.7]{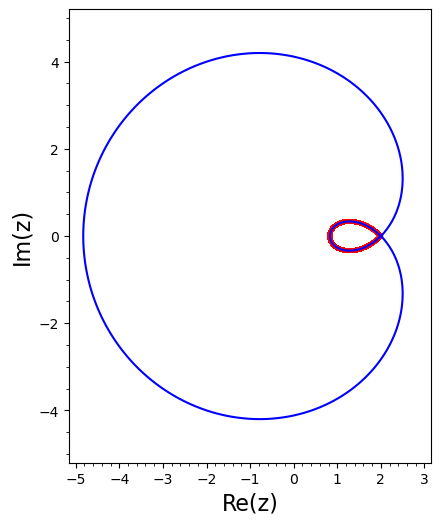}
   \caption{The set $\partial(\Omega)$ in blue and
$\partial(\sigma(c))$ in red.}
\end{figure}

\centerline{}

\end{remark}
\section{Inverse spectral mapping theorem of the quadratic Catalan equation}\label{RangeCesar}
Now we consider $(X, \Vert \quad\Vert)$ a Banach space and ${\mathcal B}(X)$ the set of linear and bounded operators. Given $T\in {\mathcal B}(X)$, as usual we write $\rho(T)$ the resolvent set given $\lambda \in \C$ such that $(\lambda-T)^{-1}\in {\mathcal B}(X)$ and $\sigma(T):=\C\backslash \rho(T)$. The spectrum of an operator is a non-empty closed set such that $\sigma(T)\subset \overline{D(0, \Vert T\Vert)}$ (\cite{[Ru]}).

In this section we study spectral properties of the solution of quadratic equation (\ref{Ceq}) with $T\in {\mathcal B}(X)$.  We say that $Y\in  {\mathcal B}(X)$ is a solution of (\ref{Ceq}) when the equality holds. Depend on $T$, the equation (\ref{Ceq}), has no solution, one, two or infinite solution, see subsection \ref{matrix}.

The proof of the following lemma is a direct consequence of the equality (\ref{Ceq}).

\begin{lemma}\label{inverse} Given $T\in {\mathcal B}(X)$ and $Y$ a solution of (\ref{Ceq}). Then $Y$ has left-inverse and $Y^{-1}_l=I-TY$.
\end{lemma}

\begin{theorem} \label{inve} Given $T\in {\mathcal B}(X)$ and $Y$ a solution of (\ref{Ceq}). Then the following are equivalent.
\begin{itemize}
\item[(i)] $0\in \rho(Y)$.
\item[(ii)] $T=Y^{-1}-Y^{-2}$.
\item[(iii)] $T$ and $ Y$ commute.
\item[(iv)] $ TY^2=YTY.$

\end{itemize}

\end{theorem}
\begin{proof} (i) As $0 \in \rho(Y)$, we obtain item (ii) from (\ref{Ceq}). The expression of $T$ in (ii) implies that $T$ and $ Y$ commute. Now, if $T$ and $ Y$ commute, then the equality $ TY^2-YTY=0$ holds. Finally, we show that item (iv) implies (i). By lemma \ref{inverse}, we have that $I-TY$ is a left-inverse of $Y$ and it is enough to check if is a right-inverse
$$
Y(I-TY)=Y-YTY=Y-TY^2=I,
$$
where we have applied (iv) and the equation (\ref{Ceq}).
\end{proof}

In the case that $dim(X)<\infty$, to be left-invertible implies to  be invertible and the conditions of Theorem \ref{inve} hold.

\begin{corollary} Let $X$ be a Banach space with $dim(X)<\infty$,  $T\in {\mathcal B}(X)$ and $Y$ a solution of (\ref{Ceq}). Then $Y$ is invertible, $T$ and $Y$ commute and $T=Y^{-1}-Y^{-2}$

\end{corollary}

In the next theorem we give the expression of $(\lambda-Y)^{-1}$ which extend the equality $Y^{-1}=I-TY$ given in Lemma \ref{inverse}.

\begin{theorem}\label{spectral} Given $T\in {\mathcal B}(X)$ and $Y$ a solution of (\ref{Ceq}) such that $0\in \rho(Y)$.
\begin{itemize}

\item[(i)] Given $\lambda \in \C$ such that $\displaystyle{\lambda-1\over \lambda^2}\in \rho(T)$ then $\lambda\in \rho(Y)$ and
$$
(\lambda-Y)^{-1}={1\over \lambda}+{1\over \lambda^3}({\lambda-1\over \lambda^2}-T)^{-1}+{Y\over \lambda^2}- {(\lambda-1)Y\over \lambda^4}({\lambda-1\over \lambda^2}-T)^{-1}.
$$

\item[(ii)] Given $\lambda \in \rho(Y)$ such that $\displaystyle{\lambda\over \lambda-1}\in \rho(Y)$ then $\displaystyle{\lambda-1\over \lambda^2}\in \rho(T)$ and
$$
({\lambda-1\over \lambda^2}-T)^{-1}={\lambda^4\over \lambda-1}({\lambda\over \lambda-1}-Y)^{-1}\left((\lambda-Y)^{-1}-{\lambda+Y\over \lambda^2}\right).
$$

\end{itemize}
\end{theorem}

\begin{proof} (i)
For $\displaystyle{\lambda-1\over \lambda^2}\in \rho(T)$, we have that
$$
{1\over \lambda}+{1\over \lambda^3}({\lambda-1\over \lambda^2}-T)^{-1}={1\over \lambda}\left(1+{1\over \lambda-1-\lambda^2T}\right)={\lambda T-1\over \lambda^2 T-\lambda+1},
$$
$$
{Y\over \lambda^2}- {(\lambda-1)Y\over \lambda^4}({\lambda-1\over \lambda^2}-T)^{-1}={Y\over \lambda^2}\left(1-{\lambda-1\over \lambda-1-\lambda^2T}\right)={YT\over \lambda^2 T-\lambda+1}.
$$
To conclude the equality, we check that
$$
(\lambda-Y)\left({\lambda T-1+YT\over \lambda^2T-\lambda+1}\right)= {\lambda^2T-\lambda+Y-TY^2\over \lambda^2 T-\lambda+1}=1,
$$
where we have applied the quadratic equation  (\ref{Ceq}).

\noindent(ii) Take  $\lambda \in \rho(Y)$ such that $\displaystyle{\lambda-1\over \lambda}\in \rho(Y)$.
$$
{\lambda^4\over \lambda-1}({\lambda\over \lambda-1}-Y)^{-1}\left((\lambda-Y)^{-1}-{\lambda+Y\over \lambda^2}\right)={\lambda^2 Y^2\over(\lambda-(\lambda-1)Y)(\lambda-Y)}.
$$
Now we check that
$$
({\lambda-1\over \lambda^2}-T){\lambda^2 Y^2\over(\lambda-(\lambda-1)Y)(\lambda-Y)}={(\lambda-1) Y^2-\lambda^2 TY^2\over(\lambda-(\lambda-1)Y)(\lambda-Y)}=1,
$$
where  we have applied  the quadratic equation  (\ref{Ceq}) in the last equality. \end{proof}

\begin{remark} The part (i) of Theorem \ref{spectral} may be considered as an inverse spectral mapping theorem: given $\lambda \in \sigma(Y)$ then $\displaystyle{\lambda-1\over \lambda^2}\in \sigma(T)$, in fact, $T=\displaystyle{Y-1\over Y^2}$, see  Theorem \ref{inve} (ii).

\end{remark}

\section{ Catalan generating functions for bounded operators}

In this section, we consider the particular case that $T$ is a linear and bounded operator on the Banach space $X$, $T\in {\mathcal B}(X)$,   such that
\begin{equation}\label{pb}
\sup_{n\ge 0} {\Vert 4^n T^n\Vert}:=M<\infty,
\end{equation}
i.e., ${\displaystyle{4T}}$ is a power-bounded operator. In this case $\sigma(T)\subset \overline{D(0,{1\over 4})}$. Under the condition (\ref{pb}), we may define the following bounded operator,
\begin{equation} \label{serie}
C(T):=\sum_{n\ge 0}C_nT^n,
\end{equation}
as a direct consequence of (\ref{norm}). Moreover, the bounded operator $C(T)$ may be consider as the  image of the Catalan sequence $c=(C_n)_{n\ge 0}$ in the algebra homomorphism $\Phi: \ell^1(\N^0, {1\over 4^n}) \to {\mathcal B}(X)$ where
$$
\Phi(a)x:=\sum_{n\ge 0} a_nT^n(x), \qquad a=(a_n)_{n\ge 0}\in  \ell^1(\N^0, {1\over 4^n}) , \quad x\in X,
$$
i.e., $\Phi(c)=C(T)$. The $\Phi$ algebra homomorphism (also called functional calculus) has been considered in several papers, two of them are \cite[Section 2]{[Dun]} and more recently \cite[Section 5.2]{[GT]}.  In particular, the map $\Phi$ allows to define the following operators
\begin{eqnarray*}
\Phi(\delta_n)&=&T^n, \qquad n\ge 0,\cr
\Phi({1\over \lambda}{ p_\lambda})&=&(\lambda-T)^{-1}, \qquad \vert \lambda\vert>{1\over 4},\cr
{\sqrt{{1}-{4T}}}&=&\sum_{n\ge 0} {(-4)^n}{ {1\over 2}\choose n} T^n,\cr
\end{eqnarray*}
where we have applied the ``generalized binomial formula'', $(1+z)^\alpha=\sum_{n\ge 0}{\alpha \choose n}z^n$, for $\vert z\vert <1$ and ${\alpha \choose n}={\alpha(\alpha-1)\cdots(\alpha-n+1)\over n!}$ for $\alpha>0$. Remind that ${\alpha \choose n}\sim {1\over n^{1+\alpha}}$ when $n\to \infty$.

\begin{theorem}\label{genne} Given  $T\in {\mathcal B}(X)$ such that ${4T}$ is power-bounded and $c=(C_n)_{n\ge 0}$ the Catalan sequence. Then
\begin{itemize}
\item[(i)] The operator $C(T)$ defined by (\ref{serie}) is well-defined, $T$ and $C(T)$ commute, and $C(T)$ is a solution of the quadratic equation (\ref{Ceq}).
\item[(ii)] The following integral representation holds
$$
C(T)x={1\over \pi}\int_{1\over 4}^\infty {\sqrt{\lambda-{1\over 4}}\over \lambda}(\lambda-T)^{-1}x d\lambda, \qquad x\in X.
$$
\item[(iii)] The following integral representation holds
$$
TC(T)={I\over 2}-{\sqrt{{1\over 4}-T}}.
$$
\item[(iv)] The spectral mapping theorem holds for $C(T)$, i.e, $\sigma(C(T))=C(\sigma(T))$ and
$$
\sigma(C(T))\subset C(\overline{D(0,{1\over 4})})\subset \sigma(c).
$$
\item[(v)] Given $\lambda \in \C$ such that $\displaystyle{\lambda-1\over \lambda^2}\in \rho(T)$ then $\lambda\in \rho(C(T))$ and
$$
(\lambda-C(T))^{-1}={1\over \lambda}+{1\over \lambda^3}({\lambda-1\over \lambda^2}-T)^{-1}+{C(T)\over \lambda^2}- {(\lambda-1)C(T)\over \lambda^4}({\lambda-1\over \lambda^2}-T)^{-1}.
$$
\end{itemize}

\end{theorem}
\begin{proof} (i) From (\ref{norm}),   $C(T)\in {\mathcal B}(X)$ as we have commented. It is clear that $T$ and $C(T)$ commute. We apply the algebra homomorphism to the  equality given in Proposition \ref{first} (iii) to get
$$
0=\Phi(\delta_1* c^{*2}-c+\delta_0)=\Phi(\delta_1)(\Phi( c))^{*2}-\Phi(c)+\Phi(\delta_0)=T(C(T))^2-C(T)+I.
$$

\noindent (ii) As the homomorphism $\Phi$ is continuous, we apply the formula (\ref{integra}) to get
$$
C(T)x= \Phi\left({1\over \pi}\int_{1\over 4}^\infty{\sqrt{\lambda -{1\over 4}}\over \lambda^{2}}p_\lambda d\lambda\right)x={1\over \pi}\int_{1\over 4}^\infty {\sqrt{\lambda-{1\over 4}}\over \lambda}(\lambda-T)^{-1}x d\lambda,
$$
for $x\in X$.

\noindent(iii) Note that
$$
{I\over 2}-{\sqrt{{1\over 4}-T}}=-{1\over 2}\sum_{n\ge 0} {(-4)^{n+1}}{ {1\over 2}\choose n+1} T^{n+1}=TC(T),
$$
since $\displaystyle{C_{n}= -{1\over 2} {(-4)^{n+1}}{ {1\over 2}\choose n+1}}$ for $n\ge 0$.

\noindent(iv) Since ${4T}$ is power bounded, the spectral mapping theorem for $C(T)$ may found in \cite[Theorem 2.1]{[Dun]}. As $\sigma(T)\subset \overline{D(0,{1\over 4})}$, we apply Proposition \ref{spec} to conclude that $\sigma(C(T))\subset C(\overline{D(0,{1\over 4})})\subset \sigma(c).$

\noindent  As $C(T)$ is a solution of  (\ref{Ceq}) such that $0\in \rho(Y)$, the item (v) is a particular case of
Theorem \ref{spectral} (i). An alternative proof may be obtained using Theorem \ref{inver} and the algebra homomorphism $\Phi$ for $\lambda \in \Omega \backslash \{0\}.$
\end{proof}

\begin{remark} In the case that $\sigma(T)\subset D(0,{1\over 4})$, the generating function $C(z)$ given in (\ref{gene}) is an holomorphic function in a neighborhood of $\sigma(T)$. Then the Dunford functional calculus, defined by the integral Cauchy-formula,
$$
f(T)x:=\int_{\Gamma}f(z)(z-T)^{-1}xdz, \qquad x\in X,
$$
 allows to defined $C(T)$, (\cite[Section VIII.7]{Y}) which, of course, coincides with the expression gives in (\ref{serie}). As usual, the path $\Gamma$ rounds the spectrum set $\sigma(T)$.
\end{remark}
\section{Examples, applications and final comments}

 In this section we present some particular examples of operators $T$  for which we solve  the equation  (\ref{Ceq}). In the subsection 6.1, we consider the Euclidean space $\C^2$ and matrices $T=\lambda I_2, \,\lambda \begin{pmatrix}
0 & 1\\
1 & 0
\end{pmatrix}, \lambda \begin{pmatrix}
0 & 1\\
0 & 0
\end{pmatrix}$ where $\lambda \in \C$. Note that we have to solve a system of four quadratic equations.  We also calculate $C(T)$ for these matrices. In  subsection 6.2  we check $C(a)$ for some $a\in \ell^1(\N^0, {1\over 4^n})$. Finally we present an iterative method for matrices $\R^{n\times n }$ which are defined using Catalan numbers in subsection 6.3.

\subsection{Matrices on $\C^2$}\label{matrix} We consider the Euclidean space ·$\C^2$ and the operator $T=\lambda I_2$ with $0\not=\lambda\in \C$. Then the solution of  (\ref{Ceq}) is given by
$$
Y=\begin{pmatrix}
{1\pm\sqrt{1-4\lambda(1+\lambda bc)}\over 2\lambda} & b\\
 c & {1\mp\sqrt{1-4\lambda(1+\lambda bc)}\over 2\lambda}
\end{pmatrix}, \qquad
$$
for  $\vert c\vert+\vert b\vert>0$ where  the allowed signs are $(-, +)$ and $(+,-)$. For $c=b=0$,    the solutions are given by
$$
Y=\begin{pmatrix}
{1\pm\sqrt{1-4\lambda}\over 2\lambda} & 0\\
0 & {1\pm\sqrt{1-4\lambda}\over 2\lambda}
\end{pmatrix}, \qquad
$$
where the allowed signs are all four combinations. In both cases, note that $\sigma(Y)=\{C(\lambda), \displaystyle{1\over \lambda C(\lambda)}\} $ and $\sigma(T)=\{\lambda\}$, compare with Theorem \ref{spectral}. In the case that $\vert \lambda \vert\le {1\over 4}$.
$$
C(T)=\begin{pmatrix}
C(\lambda) & 0\\
0 & C(\lambda)
\end{pmatrix}. \qquad
$$

Now we study the case $T=\begin{pmatrix}
0 & \lambda\\
\lambda & 0
\end{pmatrix}$ with $\lambda \in \C\backslash\{0\}$. The solutions of (\ref{Ceq}) are given by
$$
Y=\begin{pmatrix}
{a} & {a-1\over 2\lambda a}\\
{a-1\over 2\lambda a} & a
\end{pmatrix}
$$
where $a$ is a solution of the biquadratic equation $4\lambda^2a^4-a^2+1=0$. In the case that $\vert \lambda \vert\le {1\over 4}$, we get that
$$
C(T)=\begin{pmatrix}
C_e(\lambda) & C_o(\lambda)\\
 C_o(\lambda) &  C_e(\lambda)
\end{pmatrix},
$$ where functions $C_e $ and $C_o$ are defined in Proposition \ref{descopo}

Finally take now $T=\begin{pmatrix}
0 & \lambda\\
 0 & 0
\end{pmatrix}$ with $\lambda \in \C$. The only solution of (\ref{Ceq})  is given by $Y= \begin{pmatrix}
1 & \lambda\\
 0 & 1
\end{pmatrix}=C_0I_2+C_1T$; note that $T^n=0$ for $n\ge 2$.

\subsection{Catalan operators on  $\ell^p$} We consider the space of sequences $\ell^p(\N^0, {1\over 4^n})$  where
$$
\Vert a\Vert_{p, {1\over 4^n}}:=\left(\sum_{n=0}^\infty {\vert a_n\vert^p\over 4^{np}}\right)^{1\over p}<\infty,
$$
for $1\le p<\infty$ and  $\ell^\infty(\N^0, {1\over 4^n})$ the space of sequences embedded with the norm
$$
\Vert a\Vert_{\infty, {1\over 4^n}}:=\sup_{n\ge 0}{\vert a_n\vert\over 4^{n}}<\infty.
$$
 Note that $\ell^1(\N^0, {1\over 4^n})\hookrightarrow \ell^p(\N^0, {1\over 4^n})\hookrightarrow \ell^\infty(\N^0, {1\over 4^n})$.

 Now we consider $c=(C_n)_{n\ge 0}$ the Catalan sequence and the convolution operator $C(f):= c\ast f$  for $f\in \ell^p(\N^0, {1\over 4^n})$ with $1\le p\le \infty$. Since $C(f)=\sum_{n\ge 0} c_n\delta_n(f)=\sum_{n\ge 0} c_n(\delta_1)^n(f)$, we apply Theorem \ref{genne} (iv) to get
 $$
 \sigma(C)=C(\sigma(\delta_1))=C(\overline{D(0,{1\over 4})}),
 $$
 i.e., it is independent on $p$ and coincides with the spectrum of the Catalan sequence $c$ in  $\ell^1(\N^0, {1\over 4^n})$
 (Proposition \ref{spec}).

 Now we consider the spaces $\ell^p(\Z)$ for $1\le p\le \infty$ defined in the usual way. The element $a=\delta_1-\delta_0$ defines the classical backward difference operator
 $$a\ast (f)(n):= f(n-1)-f(n), \qquad f\in \ell^p(\Z),
 $$
 for $n\in \Z$. Note that $\Vert a\Vert =2$, and
 $$
 (\lambda\delta_0+a)^{-1}=\sum_{j\ge 0}{\delta_j\over (1+\lambda)^{j+1}}, \quad 1<\vert 1+\lambda\vert,
 $$
 see \cite[Theorem 3.3 (4)]{[GLM]}. Now we need to consider $\displaystyle{a\over 8}$ and the associated Catalan generating operator defined by (\ref{serie}). By Theorem \ref{genne} (2), we get that
 \begin{eqnarray*}
 C({a\over 8})&=&{8\over \pi}\int_{1\over 4}^\infty {\sqrt{\lambda-{1\over 4}}\over \lambda}(8\lambda\delta_0+a)^{-1}x d\lambda\cr&=& {4\over \pi}\sum_{j\ge 0}{\delta_j\over 2^{j+1}}\int_0^\infty{\sqrt{t}\over (t+1)(t+{3\over 2})^{j+1}}dt\cr
 &=&(2\sqrt{6}-4)\delta_0+\sum_{j=1}^\infty\left({\sqrt{6}\over 3}\sum_{k=j}^{\infty}{C_k\over 12^k}\right)\delta_j,
 \end{eqnarray*}
where we have applied Theorem \ref{integrals} for $z={3\over 2}$.

A similar results holds for the forward difference operator defined by $\Delta f (n):=f(n+1)-f(n)$, \cite[Theorem 3.2]{[GLM]}.

\medskip

\subsection{Iterative methods  on $\R^n$ applied to Quasi-birth-death processes}

The quadratic matrix equation:
 \begin{equation}\label{QME}
\mathcal{Q}(Y):=  TY^2- Y+I=0,
\end{equation}
  is related to the particular Markov chain   characterized  by   its transition matrix $P$ which is an infinite block  tridiagonal matrix of the form:
$$P=  \left(\begin{array}{ccccc }
D_1& D_2& & &0 \\
I& 0&T&& \\
 &I& 0&T&\\
0&& \ddots&\ddots&\ddots
\end{array}\right),$$
  where the blocks $D_1$, $D_2$, $I$, $T$ are $n\times n$ nonnegative matrices such that  $D_1+D_2$ and  $I+T$ are row stochastic.   A discrete-time Markov chain represent
a quasi-birth-death stochastic process. In fact  a quasi-birth-death stochastic process is a Discrete-time Markov chain having infinitely states (\cite{LAT}). Thus, a nonnegative solution of quadratic matrix equation
 (\ref{QME}) is necessary to describe probabilistically the behavior of that Markov chain.

In \cite{D1, D2} the author demonstrated the usefulness of Newton's method for solving the quadratic matrix equation. There are many papers containing algorithmic methodologies and
acceleration techniques related to quadratic matrix equations, see for instance \cite{D1, D2, HRR, LAN85, Rogers}.

Our purpose in this section
is to show experimentally the benefits of  a higher order iterative method to approximate the nonnegative solution of  equation
 (\ref{QME}) which uses the Catalan numbers, $C_j$:
\begin{equation}\label{Cata}
 \left\{ \begin{array}{ll}
Y_0\ \text{given},\\[1ex]
Y_{n+1}= Y_n +\displaystyle\sum_{j\ge 0}H_j(Y_n),
 \quad n\geq 0,
\end{array} \right.\end{equation}
with $$H_k(Y)=-\frac{C_k}{2^k}\left([\mathcal{Q}'(Y)]^{-1}\mathcal{Q}''(Y)[\mathcal{Q}'(Y)]^{-1}\mathcal{Q}(Y) \right)^k[\mathcal{Q}'(Y)]^{-1}\mathcal{Q}(Y). $$
Notice that method (\ref{Cata}) has infinite speed of convergence to approximate  a solution of equation (\ref{QME}), see \cite{HR}. That is, the  solution is obtained in the first iteration. To apply this method carries on computing the square root of the matrix $(I-4T)$. To avoid this, we can
truncate the series, thus obtaining a high-order method of convergence.

This method can be write in terms of Sylvester equations, one of the most often in matrix equations  (\cite{LAN85}):
 \begin{equation}\label{Lyapunov}
AY + YB =D,
\end{equation}
with  $A$, $B$ and  $D$ given matrices.

Taking into account that the first  Fr\'{e}chet  derivative at a matrix $Y$  is a linear map $\mathcal{Q}'(Y): \R^{m\times m}   \rightarrow \R^{m\times m}$ such that
\begin{equation}\label{primeraderivada}
\mathcal{Q}'(Y) E=
 TEY+(TY- I) E,
\end{equation}
and  the second  derivative at $Y$,  $ \mathcal{Q}''(Y):  \R^{m\times m}\times  \R^{m\times m}    \rightarrow \R^{m\times m}  $ is given by
\begin{equation}\label{segderivada}
\mathcal{Q}''(Y) E_1 E_2=
 T(E_1E_2+  E_2E_1),
\end{equation}
 is a bilinear constant operator.

Notice that  method (\ref{Cata})   can be written in terms of generalized Sylvester equations, for $j\ge 1$,
\begin{equation}\label{Catalan}
 \left\{
 \begin{array}{lll}
Y_0\ \text{given},\\[1ex]
 ( TY_n-I)H_0 +  H_0 T Y_n  &=&  -TY_n ^2-  (Y_n-I), \\[1ex]
 (T Y_n-I)H_j +  H_j TY_n  &=& -\displaystyle\frac{C_j}{2C_{j-1}}T(H_0 H_{j-1} +H_{j-1} H_0 ),\\[1ex]
  (TY_n-I) Y_{n+1}+ Y_{n+1}TY_n  &=& \\[1ex] \hfil \qquad \qquad \qquad\qquad \qquad TY_n^2-I&-&\displaystyle{\sum_{j\ge 1}\frac{C_j}{2C_{j-1}}T(H_0 H_{j-1} +H_{j-1} H_0 )}. \\[1ex]
  \end{array}
 \right.
\end{equation}

Notice that, method (\ref{Catalan4}) is reduced to solve   three Sylvester equations with the same matrix system.
The Bartels-Stewart algorithm is ideally suited to the sequential solution of Sylvester equation (\ref{Lyapunov}) with the same matrix system (\cite{Bartels}).

In particular, method of (\ref{Cata})  truncated to $k=2$ has   forth order of converge:
\begin{equation}\label{Catalan4}
 \left\{
 \begin{array}{lll}
Y_0\ \text{given},\\[1ex]
( TY_n-I)H_0 +  H_0 T Y_n   &=&  -T Y_n^2-  (Y_n-I), \\[1ex]
(TY_n-I)(H_0+H_1)+(H_0+H_1)T Y_n    &=& -TY_n ^2-  (Y_n-I)-T H_0^2,\\[1ex]
(TY_n-I) Y_{n+1}+ Y_{n+1}T Y_n  &=& TY_n^2-I-\frac{1}{2}T \left(H_0(H_0+H_1)\right.
\\[1ex]&&\left.+(H_0+H_1)H_0\right).
 \end{array}
 \right.
\end{equation}

Taking into account that the most commonly used  Newton's method:
$$(TY_n-I) Y_{n+1}+ Y_{n+1}TY_n  = TY_n^2,
$$
only achieve quadratic convergence speed   (\cite{D1, D2}), the forth order method (\ref{Catalan4})    is a good alternative to approximate a solution of quadratic equation (\ref{QME}).

Next,  a numerical example is shown where the matrix $T$ is ill conditioned. With high accuracy we approximate numerically  the nonnegative solution    of equation (\ref{QME}) using the   method (\ref{Catalan4}). To do that, we take  $T=(t_{ii})$ a diagonal matrix $100\times 100$ with entries $t_{ii}= 10^{-1}, i=1\ldots, 9$ and  $t_{10,10}= 10^{-10}$.
Method (\ref{Catalan4}) is implemented   in Mathematica Version $10.0$, with stopping criterion $RES <10^{-20}$, $RES :=\|\   \mathcal{Q}(Y_k)\|_\infty$.
We choose   the starting matrix $Y_0= T.$ We show the number of iterations necessary to achieve the required precision.
The numerical results are reported in Table  \ref{res}.

 \begin{table}[hhh!!!]
\begin{tabular}{ p{2cm} p{4cm}p{4cm}p{4cm} }
\hline\noalign{\smallskip}
  $k$ &Newton & method  (\ref{Catalan4})\\
 \hline
$1 $  & $8.45274\ldots \times 10^{-2}$  &$  1.03079\ldots  \times 10^{-2 }  $    \\
$2$   &$1.12729\ldots \times 10^{-3}$  &$3.01635\ldots \times 10^{-8} $     \\
$3$  &$2.11638\ldots \times 10^{-7}$  &$7.62333\ldots \times 10^{-25}$     \\
$4$  &$7.46507\ldots \times 10^{-15}$  &    \\
$5$  &$9.28789\ldots \times 10^{-30}$  &     \\
  \hline
\noalign{\smallskip}\hline\noalign{\smallskip}
\end{tabular}
\caption{Residuals  for  the Newton method and  method  (\ref{Catalan4}), with stopping criteria $\|\ \mathcal{Q}(Y_k)\|_\infty< 10^{-20}$.
\label{res} }
\end{table}



\end{document}